\documentclass[12pt]{article}
\usepackage[a4paper]{geometry}
\usepackage{amssymb,amsmath,amsthm}
\usepackage{graphicx,cite,xcolor}
\usepackage{longtable,pdflscape,booktabs,caption,multicol}
\usepackage[colorlinks=true,citecolor=black,linkcolor=black,urlcolor=blue]{hyperref}
\usepackage{enumitem}

\theoremstyle{plain}
\newtheorem{theorem}{Theorem}
\newtheorem{lemma}[theorem]{Lemma}
\newtheorem{corollary}[theorem]{Corollary}
\newtheorem{question}[theorem]{Question}

\newtheorem{conjecture}[theorem]{Conjecture}

\theoremstyle{definition}

\theoremstyle{remark}

\newcommand{\arxiv}[1]{\href{https://arxiv.org/abs/#1}{\texttt{arXiv:#1}}}
\newcommand{\doi}[1]{\href{https://doi.org/#1}{\texttt{doi:#1}}}

\usepackage{listings}
\definecolor{mauve}{rgb}{0.58,0,0.82}
\lstdefinestyle{pitonche} {
    language = Python,
    basicstyle = footnotesizettfamily,
    showspaces = false,
    showstringspaces = false,
    breakautoindent = true,
    flexiblecolumns = true,
    keepspaces = true,
    stepnumber = 1,
    xleftmargin = 0pt
}
\usepackage{float}
\usepackage{tikz}
\usetikzlibrary{arrows}
\usetikzlibrary{arrows.meta}
\usetikzlibrary{chains}
\usetikzlibrary{positioning}
\usetikzlibrary{automata,positioning,calc}
\usetikzlibrary{decorations}
\usetikzlibrary{decorations.shapes}
\usetikzlibrary{decorations.markings}
\tikzset{
    edge/.style={-{Latex[scale=1.7]}},
    dedge/.style={{Latex[scale=1.7]}-{Latex[scale=1.7]}},
}
\usepackage{bm}

\newcommand\Osq{\mathbin{\text{\scalebox{.84}{$\square$}}}}

\title{A note on Cayley nut graphs\\ whose degree is divisible by four}

\author{Ivan Damnjanovi\'c\thanks{The author is supported by Diffine LLC.}\\
\small University of Ni\v s, Faculty of Electronic Engineering,\\[-0.4ex]
\small Aleksandra Medvedeva 14, 18106 Ni\v s, Serbia\\[-0.4ex]
\small\tt ivan.damnjanovic@elfak.ni.ac.rs\\
\small Diffine LLC\\[-0.4ex]
\small 3681 Villa Terrace, San Diego, CA 92104, USA\\[-0.4ex]
\small\tt ivan@diffine.com}

\begin{document}

\maketitle

\begin{abstract}
A nut graph is a non-trivial simple graph such that its adjacency matrix has a one-dimensional null space spanned by a full vector. It was recently shown by the authors that there exists a $d$-regular circulant nut graph of order $n$ if and only if $4 \mid d, \, 2 \mid n, \, d > 0$, together with $n \ge d + 4$ if $d \equiv_8 4$ and $n \ge d + 6$ if $8 \mid d$, as well as $(n, d) \neq (16, 8)$ [\arxiv{2212.03026}, 2022]. In this paper, we demonstrate the existence of a $d$-regular Cayley nut graph of order $n$ for each $4 \mid d, \, d > 0$ and $2 \mid n, \, n \ge d + 4$, thereby resolving the existence problem for Cayley nut graphs and vertex-transitive nut graphs whose degree is divisible by four.

\bigskip\noindent
{\bf Mathematics Subject Classification:} 05C50, 05C25.\\
{\bf Keywords:} Cayley graph, vertex-transitive graph, circulant graph, nut graph, graph spectrum, graph eigenvalue.
\end{abstract}

\newpage
\section{Introduction}

A \emph{nut graph} represents a non-trivial simple graph whose adjacency matrix contains a one-dimensional null space spanned by a full vector. These graphs were introduced by Sciriha \cite{Sciriha1, Sciriha2, Sciriha3, Sciriha4} and their properties were subsequently investigated in \cite{Sciriha5, Sciriha6}. Furthermore, the chemical justification for studying such graphs can be found in \cite{Sciriha7, Fowler1, Coolsaet, Fowler2}. Besides that, many other results concerning nut graphs are to be found in the monograph \cite{ScirihaFarrugia}.

The study of regular nut graphs was initiated by Sciriha and Fowler in \cite{Sciriha7}. Afterwards, Gauci et al.\ \cite{Gauci} determined all the orders that a cubic or quartic nut graph could have --- there exists a cubic nut graph of order $n$ if and only if $n = 12$ or $2 \mid n, \, n \ge 18$, and there exists a quartic nut graph of order $n$ if and only if $n \in \{8, 10, 12\}$ or $n \ge 14$. The said result was later extended by Fowler et al.\ \cite{Fowler3}, who found all the orders that a $d$-regular nut graph can have for each $d \le 11$. In the same paper, Fowler et al.\ initiated the study of vertex-transitive nut graphs by asking the following question.

\begin{question}[{Fowler et al.\ \cite[Question 4]{Fowler3}}]\label{vertex_transitive_nut_question}
    For what pairs $(n, d)$ does a vertex-transitive nut graph of order $n$ and degree $d$ exist?
\end{question}

Alongside Question \ref{vertex_transitive_nut_question}, some of the corresponding necessary conditions were given in the form of the following theorem.

\begin{theorem}[{Fowler et al.\ \cite[Theorem 5]{Fowler3}}]\label{fowler_th}
    Let $G$ be a vertex-transitive nut graph on $n$ vertices, of degree $d$. Then $n$ and $d$ satisfy the following conditions --- either $d \equiv_4 0$ and $n \equiv_2 0$ and $n \ge d + 4$; or $d \equiv_4 2$ and $n \equiv_4 0$ and $n \ge d + 6$.
\end{theorem}

Subsequently, Bašič et al.\ \cite{Basic} demonstrated that there exists a $12$-regular nut graph of order $n$ if and only if $n \ge 16$ and disclosed the next two conjectures regarding the circulant nut graphs.

\begin{conjecture}[{Bašić et al. \cite[Conjecture 3.2]{Basic}}]\label{basic_conjecture_1}
    For every $d$, where $d \equiv_4 0$, and for every even $n, \, n \ge d + 4$, there exists a circulant nut graph of order $n$ and degree $d$.
\end{conjecture}
\begin{conjecture}[{Bašić et al. \cite[Conjecture 3.3]{Basic}}]\label{basic_conjecture_2}
    For every even $n, \, n \ge 16$, there exists a circulant nut graph of order $n$ and degree $12$.
\end{conjecture}

In a series of papers \cite{Damnjanovic1, Damnjanovic2, Damnjanovic3}, the authors resolved Conjectures \ref{basic_conjecture_1} and \ref{basic_conjecture_2} by finding all the pairs $(n, d)$ for which there exists a $d$-regular circulant nut graph of order $n$, as shown in the following theorem.

\begin{theorem}[{Damnjanović \cite[Theorem 5]{Damnjanovic3}}]\label{circulant_nut_resolution}
    For each $d \in \mathbb{N}_0$ and $n \in \mathbb{N}$, there exists a $d$-regular circulant nut graph of order $n$ if and only if the next three conditions are satisfied:
    \begin{itemize}
        \item $ d > 0, \, 4 \mid d, \, 2 \mid n$;
        \item $n \ge d + 4$ if $d \equiv_8 4$, and $n \ge d + 6$ if $8 \mid d$;
        \item $(n, d) \neq (16, 8)$.
    \end{itemize}
\end{theorem}

Here, we heavily rely on Theorem \ref{circulant_nut_resolution} in order to investigate the structural properties of Cayley nut graphs and provide a complete answer to Question \ref{vertex_transitive_nut_question} for the case $4 \mid d$. The main result of the paper is given in the next theorem.

\begin{theorem}\label{main_th}
    For each $d \in \mathbb{N}, \, 4 \mid d$, there exists a $d$-regular Cayley nut graph of each even order $n$ such that $n \ge d + 4$.
\end{theorem}

The rest of the paper will primarily focus on providing a full mathematical proof of Theorem \ref{main_th} and discussing some of its immediate corollaries. All graphs considered will be simple and finite and standard notation shall be used. Also, all spectral properties considered will correspond to the standard $(0, 1)$-adjacency matrix of a simple graph.

\section{Main result}

For starters, given the fact that each circulant graph of order $n$ represents a Cayley graph corresponding to the cyclic group $\mathbb{Z}_n$, it is clear that whenever there exists a $d$-regular circulant nut graph of order $n$, there must also exist a $d$-regular Cayley nut graph of order $n$. Bearing in mind Theorem \ref{circulant_nut_resolution}, it becomes evident that in order to finalize the proof of Theorem \ref{main_th}, it is sufficient to demonstrate the existence of a $d$-regular Cayley nut graph of order $n$ when:
\begin{enumerate}[label={\bf (\roman*)}]
    \item\label{case_1} $8 \mid d$ and $n = d + 4$; or
    \item $(n, d) = (16, 8)$.
\end{enumerate}
We will cover case \ref{case_1} in the next lemma.

\begin{lemma}
    For each $d \in \mathbb{N}, \, 8 \mid d$, the graph $\overline{C_{\frac{d+4}{2}} \Osq K_2}$ represents a $d$-regular Cayley nut graph of order $d + 4$.
\end{lemma}
\begin{proof}
    Without loss of generality, let $\beta = \frac{d+4}{2}$ and let the vertices of $G = C_{\frac{d+4}{2}} \Osq K_2$ be denoted via $a_j, b_j, j = \overline{0, \beta-1}$, so that the neighbors of each $a_j$ are $a_{j-1}, a_{j+1}, b_j$ and the neighbors of each $b_j$ are $b_{j-1}, b_{j+1}, a_j$, where all the index arithmetics are performed modulo $\beta$. It is clear that $G$ represents a cubic Cayley graph corresponding to the group $\mathbb{Z}_\beta \times \mathbb{Z}_2$. From here, it immediately follows that $\overline{G}$ is a $d$-regular Cayley graph of order $d + 4$. Thus, in order to complete the proof of the lemma, it is enough to show that $\overline{G}$ is of nullity one, since each Cayley graph is necessarily vertex-transitive, and it is easy to check that vertex-transitive graphs of nullity one must be nut.

    Let $\bm{u} \in \mathbb{R}^{V_G}$ be a vector from the null space of $\overline{G}$, i.e.\ a vector such that
    \begin{equation}\label{aux_1}
        \bm{A}_{\overline{G}} \, \bm{u} = \bm{0} .
    \end{equation}
    It is obvious that $\bm{u} \perp \bm{1}$, where $\bm{1}$ denotes the all ones vector. For this reason, we have
    \[
        \sum_{j = 0}^{\beta} \bm{u}_{a_j} + \sum_{j = 0}^{\beta} \bm{u}_{b_j} = 0,
    \]
    hence Eq.\ (\ref{aux_1}) promptly gives us
    \begin{alignat}{2}
        \label{aux_2} \bm{u}_{a_j} + \bm{u}_{a_{j-1}} + \bm{u}_{a_{j+1}} + \bm{u}_{b_j} &= 0 \qquad \left( \forall j = \overline{0, \beta-1} \right),\\
        \label{aux_3} \bm{u}_{b_j} + \bm{u}_{b_{j-1}} + \bm{u}_{b_{j+1}} + \bm{u}_{a_j} &= 0 \qquad \left( \forall j = \overline{0, \beta-1} \right) .
    \end{alignat}
    By subtracting Eq.\ (\ref{aux_3}) from Eq.\ (\ref{aux_2}), we reach
    \begin{alignat*}{2}
        && \bm{u}_{a_{j-1}} + \bm{u}_{a_{j+1}} - \bm{u}_{b_{j-1}} - \bm{u}_{b_{j+1}} &= 0 \\
        \implies \quad && \bm{u}_{a_{j+1}} - \bm{u}_{b_{j+1}} &= -\left( \bm{u}_{a_{j-1}} - \bm{u}_{b_{j-1}} \right) ,
    \end{alignat*}
    which further leads us to
    \begin{equation}\label{aux_4}
        \bm{u}_{a_{2j}} - \bm{u}_{b_{2j}} = (-1)^j \left( \bm{u}_{a_0} - \bm{u}_{b_0} \right)
    \end{equation}
     for each $j \in \mathbb{Z}$. Since $2 \nmid \frac{d+4}{4}$, plugging in $j = \frac{d+4}{4} = \frac{\beta}{2}$ to Eq.\ (\ref{aux_4}) helps us obtain
     \begin{alignat*}{2}
        && \bm{u}_{a_0} - \bm{u}_{b_0} &= - \left( \bm{u}_{a_0} - \bm{u}_{b_0} \right)\\
        \implies \quad && \bm{u}_{a_0} - \bm{u}_{b_0} &= 0,
     \end{alignat*}
     which means that $\bm{u}_{a_0} = \bm{u}_{b_0}$. Analogously, it can be shown that
     \begin{equation}\label{aux_5}
        \bm{u}_{a_j} = \bm{u}_{b_j} \qquad \left( j = \overline{0, \beta-1} \right) .
    \end{equation}
    
    Bearing in mind Eq.\ (\ref{aux_5}), it is trivial to notice that Eqs.\ (\ref{aux_2}) and (\ref{aux_3}) simultaneously get down to
    \begin{equation}\label{aux_6}
        \bm{u}_{a_{j-1}} + \bm{u}_{a_{j+1}} = -2 \bm{u}_{a_j} \qquad \left( \forall j = \overline{0, \beta-1} \right) .
    \end{equation}
    If we denote $V_1 = \{ a_0, a_1, \ldots, a_{\beta-1} \}$, then it is not difficult to observe from Eq.~(\ref{aux_6}) that the vector $\bm{u}|_{V_1}$ can be thought of as an element of a $C_\beta$ eigenspace corresponding to the eigenvalue $-2$ (see, for example, \cite[Section 1.4.3]{Brouwer}). Since $2 \mid \beta$, it swiftly follows that
    \begin{equation}\label{aux_7}
        \bm{u}_{a_j} = (-1)^j \, \bm{u}_{a_0} \qquad \left( j = \overline{0, \beta-1} \right) .
    \end{equation}
    It is trivial to check that all vectors $\bm{u} \in \mathbb{R}^{V_G}$ satisfying the conditions (\ref{aux_5}) and (\ref{aux_7}) truly belong to the null space of $\overline{G}$. However, this set of vectors clearly forms a one-dimensional vector subspace, hence the graph $\overline{G}$ must be of nullity one.
\end{proof}

In order to finalize the proof of Theorem \ref{main_th}, it is enough to demonstrate the existence of an $8$-regular Cayley nut graph of order $16$. If we use
\[
    \mathrm{QD}_{16} = \langle s, t \colon s^8 = t^2 = 1, \, tst = s^3 \rangle
\]
to denote the quasihedral group of order $16$ as done so in \cite{quasihedral}, then the Cayley graph corresponding to this group, with the generator set
\[
    \{ tst, tsts^2, t, s^2 t, ts, sts, ts^2, sts^2 \} ,
\]
is surely $8$-regular and of nullity one. This can easily be checked via computer with the help of any mathematical software. From here, it can directly be established that an $8$-regular Cayley nut graph of order $16$ exists.

\section{Conclusion}

Taking into consideration Theorem \ref{fowler_th}, it is clear that Theorem \ref{main_th} leads us to the following two corollaries.

\begin{corollary}\label{Cayley_cor}
    For each $d \in \mathbb{N}, \, 4 \mid d$, there exists a $d$-regular Cayley nut graph of order $n$ if and only if $2 \mid n$ and $n \ge d + 4$.
\end{corollary}
\begin{corollary}\label{vt_cor}
    For each $d \in \mathbb{N}, \, 4 \mid d$, there exists a $d$-regular vertex-transitive nut graph of order $n$ if and only if $2 \mid n$ and $n \ge d + 4$.
\end{corollary}

Corollaries \ref{Cayley_cor} and \ref{vt_cor} dictate the sets of all the orders that a $d$-regular Cayley nut graph and $d$-regular vertex-transitive nut graph can have, whenever $4 \mid d$. In particular, Corollary~\ref{vt_cor} provides a partial answer to Question \ref{vertex_transitive_nut_question} by covering the case $4 \mid d$. Besides that, it is trivial to notice that a $2$-regular nut graph does not exist, given the fact that such a graph is necessarily composed of cycles. Hence, we end the paper with the next two open problems.

\begin{question}
    For each $d \in \mathbb{N}, \, d \equiv_4 2, \, d \ge 6$, find all the integers $n \in \mathbb{N}$ for which there exists a $d$-regular Cayley nut graph of order $n$.
\end{question}
\begin{question}
    For each $d \in \mathbb{N}, \, d \equiv_4 2, \, d \ge 6$, find all the integers $n \in \mathbb{N}$ for which there exists a $d$-regular vertex-transitive nut graph of order $n$.
\end{question}

\end{document}